\documentclass{article}
\usepackage{amsmath,amstext, amsthm, amssymb}
\numberwithin{equation}{section}

\DeclareMathOperator{\supp}{supp}

\def\<{\langle}
\def\>{\rangle}

\def\<{\langle}
\def\>{\rangle}

\newcommand{\sR}{{\mathbb R}}

\newcommand{\be}{\begin{equation}}
\newcommand{\ee}{\end{equation}}

      \newtheorem{theorem}{Theorem}[section]
       \newtheorem{proposition}[theorem]{Proposition}
       \newtheorem{corollary}[theorem]{Corollary}
       \newtheorem{lemma}[theorem]{Lemma}
\newtheorem{claim}[theorem]{Claim}
\theoremstyle{remark}
       \newtheorem{remark}{Remark}[section]
\theoremstyle{definition}

\def\<{\langle}
\def\>{\rangle}

\def\<{\langle}
\def\>{\rangle}

\newtheorem{CT}{Theorem}

\title{Projection formulas for orthogonal polynomials}
\author{W{\l}odzimierz Bryc\thanks{\noindent Research partially supported by NSF
grants \#INT-03-32062, \#DMS-05-04198, and by the C.P. Taft Memorial
Fund.}\\ Department of Mathematical Sciences \\
University of Cincinnati\\
Cincinnati, OH 45221-0025 \and Wojciech Matysiak\\Faculty of
Mathematics and Information Science\\
 Warsaw University of
Technology\\ pl. Politechniki 1\\
 00-661 Warszawa, Poland \and Ryszard Szwarc\thanks{Research
partially supported by KBN (Poland) under grant 2 P03A 028 25}
\\ Institute of Mathematics\\
University of Wroc\l aw,\\ pl.\ Grunwaldzki 2/4, 50-384 Wroc\l aw, Poland
\\  and
\\ Institute of Mathematics and Computer Science,\\ University of Opole, ul. Oleska 48,
45-052 Opole, Poland
\and Jacek Weso{\l}owski
 \\ Faculty of Mathematics and Information Science\\
Warsaw University of Technology\\ pl. Politechniki 1\\ 00-661
Warszawa, Poland}
\date{
June 3, 2006; revised \today}
\begin{document}

\maketitle



\begin{abstract}
We prove  a projection formula for the four-parameter family of
orthogonal polynomials that are a reparameterization of the polynomials
in the Askey-Wilson class.  By
carefully analyzing  the recurrence relations we manage to avoid
using the explicit expression for the orthogonality measure,
which would be cumbersome due to
the complexity of the reparameterization.
\end{abstract}

\section{Introduction} Projection formulas of the type
\begin{equation}\label{pr}
q_n(x)=\int p_n(y) \nu_x(dy),
\end{equation}
where $\{\nu_x\}$ is a family of probability measures, are of interest in the theory
of orthogonal polynomials and in probability.

Explicit formulas for the measure $\nu_x$ have been known since
\cite{Askey-Fitch69} when $q_n(x)$ and  $p_n(y)$ are both Jacobi
polynomials. These formulas were extended to pairs of Askey-Wilson
polynomials in \cite{Rahman85,Rahman88} and to pairs of associated
Askey-Wilson polynomials in \cite{Rahman97}. The proofs rely on
explicit evaluation of certain integrals, which is a topic of
independent interest.

Projection formulas of the type \eqref{pr} were used as a basis of
construction of certain Markov processes in
\cite{Bryc-Wesolowski-03,Bryc-Wesolowski-04,
Bryc-Matysiak-Wesolowski-04b,Bryc-Matysiak-Wesolowski-05}. The
technique of proof in these papers is less constructive and relies
on an implicit definition of 
the probability measure $\nu_x$ as the
orthogonality measure of the auxiliary family of orthogonal
polynomials. With the exception of
\cite{Bryc-Matysiak-Wesolowski-05}, these projection formulas dealt
with the pairs of polynomials within the Askey-Wilson class and in
fact differ from \cite{Rahman85,Rahman88} only in the allowed ranges
for the parameters. The purpose of this note is to provide a related
projection formula that covers one more parameter,
but also falls into the Askey-Wilson class. Our method
does not rely on the knowledge of explicit orthogonality measures
and has a more combinatorial character.

\tolerance=2000 Our goal is to analyze in detail the family of
orthogonal polynomials
$\overline{p}_n(y;t)=\overline{p}_n^{(\eta,\theta,\tau,q)}(y;t)$
which appeared   in the study of stochastic processes with linear
regressions and quadratic conditional variances in \cite[Theorem
4.5]{Bryc-Matysiak-Wesolowski-04}.  Let $\overline{p}_{-1}=0$,
$\overline{p}_0=1$. Fix $\eta,\theta\in \sR$, $\tau\geq 0$,
$-1<q\leq 1$. For $t>0$, $n\geq 0$ let
\begin{equation}
\label{p-monic} y\overline{p}_n(y;t)= \overline{p}_{n+1}(y;t)+b_n(t)
\overline{p}_n(y;t)+a_{n-1}c_n(t) \overline{p}_{n-1}(y;t),
\end{equation}
where for $\eta\ne 0$
\begin{eqnarray}
a_n&=&\eta^{-1}+\theta [n]_q+[n]_q^2\eta\tau,\label{aa}
\\
b_n(t)&=&\left(t \eta+\theta + ([n]_q+[n-1]_q) \eta\tau\right)[n]_q,
\label{bb}
\\
c_n(t)&=&\eta( t +\tau [n-1]_q)[n]_q\;.\label{cc}
\end{eqnarray}
For $\eta=0$ we need to interpret $a_{n-1}c_n(t)$  as $( t +\tau
[n-1]_q)[n]_q$. Our reason for the separation of $\eta\eta^{-1}$ between two factors
is that for $\eta>0$ we have
\begin{equation}
  \label{sum=1/eta}
  b_n(t)=a_n+c_n(t)-\frac{1}{\eta},
\end{equation} a property which will be
exploited later on.
 We use the   notation
\begin{eqnarray*}
{[n]_{q}} =1+q+\dots +q^{n-1}, \;{[n]_{q}!} &=&[1]_{q}[2]_{q}\dots
[n]_{q},
 \left[ \begin{array}{c}n \\ k \end{array}\right]_{q} =\frac{[n]_{q}!}{[n-k]_{q}![k]_{q}!},
\end{eqnarray*}%
with the usual conventions $[0]_{q}=0,[0]_{q}!=1$.

Throughout this paper, by $\mu_t$ we denote the orthogonality
measure of polynomials $\{\overline{p}_n(y;t)\}$. A sufficient
condition for existence of such a probability measure is that
$\eta\theta\geq 0$, $\tau\geq 0$, and $0\leq q\leq 1$. It is
plausible that our results are valid for a more general range of the
parameters (compare \cite{Bryc-Matysiak-Wesolowski-04b} and
\cite{Bryc-Matysiak-Wesolowski-05}), but an attempt to cover such a
range is likely to lead to additional technical complications which
should be avoided in a paper that already has a significant degree
of computational complexity.

We now compare the polynomials defined by \eqref{p-monic} with the
monic Askey-Wilson \cite{Askey-Wilson-85} polynomials $\overline{w}_n$ for
the ``generic" values of parameters
of the recurrences. Recall that polynomials $\overline{w}_n$
 are defined by the recurrence
\begin{equation}
  \label{AW}
(x-\frac12(a+a^{-1}))
\overline{w}_n(x)= \overline{w}_{n+1}(x)- \frac12(A_n+C_n)\overline{w}_n+
\frac14 A_{n-1}C_n \overline{w}_{n-1}(x),\; n\geq 0,
\end{equation} where
$$
A_n=\frac{(1-abcd q^{n-1})(1-ab q^n)(1-ac q^n)(1-ad q^n)}{a(1-abcd
q^{2n-1})(1-abcd q^{2n})},
$$
$$
C_n=\frac{a(1-q^n) (1-bc q^{n-1})(1-bd q^{n-1})(1-cd q^{n-1})}{
(1-abcd q^{2n-2})(1-abcd q^{2n-1})}.
$$
This form of the Askey-Wilson recurrence is a minor rewrite of the recurrence
in \cite[(4.3)]{Ismail-Mason-95}. The
initial condition is the usual $\overline{w}_{-1}=0$ and $\overline{w}_0=1$.

If we multiply \eqref{AW} by $ (-2\alpha)^{n+1}$, substitute
$y=-2\alpha( x-(a+a^{-1})/2)$,  and introduce  polynomials
 $\overline{p}_n(y):=(-2\alpha)^n\overline{w}_n(x)$,
we get the following recurrence
\begin{equation}
  \label{AW1} y\overline{p}_n(y)=\overline{p}_{n+1}(y)+\alpha(A_n+C_n)\overline{p}_n(y)
+\alpha^2 C_n A_{n-1}\overline{p}_{n-1}(y).
\end{equation}

On the other hand,
\begin{multline*}
  a_n=\frac{(1-q)^2+\eta\theta(1-q)+\eta^2\tau}{\eta(1-q)^2}\Big(
  1-\frac{\eta\theta(1-q)
  +2\eta^2\tau}{(1-q)^2+\eta\theta(1-q)+\eta^2\tau}q^n\\+\frac{\eta^2\tau}{(1-q)^2
  +\eta\theta(1-q)+\eta^2\tau} q^{2n}\Big),
\end{multline*}
and
\begin{equation*}
  c_n(t)= \eta\frac{(1-q)t+\tau}{(1-q)^2}\left(1-\frac{\tau}{(1-q)t+\tau}q^{n-1}\right)(1-q^n).
\end{equation*}
This can be written as $a_n=\alpha_AA_n$, and $c_n(t)=\alpha_CC_n$,
where the Askey-Wilson parameters are $d=0$, and the remaining
 parameters $a,b,c$ are determined from the system of three
equations
\begin{equation}
  \label{EqA}
  a(b+c)=\frac{\eta\theta(1-q)
  +2\eta^2\tau}{(1-q)^2+\eta\theta(1-q)+\eta^2\tau}
\end{equation}
\begin{equation}
  \label{EqC}
  a^2bc=\frac{\eta^2\tau}{(1-q)^2
  +\eta\theta(1-q)+\eta^2\tau},\;
bc=\frac{\tau}{(1-q)t+\tau}.
\end{equation}
The auxiliary coefficients are then
\begin{equation}
  \label{EqAlpha}
\alpha_A=\frac{(1-q)^2+\eta\theta(1-q)+\eta^2\tau}{\eta(1-q)^2}a,\;
\alpha_C=\eta\frac{(1-q)t+\tau}{(1-q)^2 a}.
\end{equation}
From this and \eqref{sum=1/eta} we see that
 recurrence \eqref{p-monic} can be  re-written using
a new  variable $x=y+1/\eta$ and polynomials $\bar{r}_n(x):=\bar{p}_n(y)$
as
\begin{equation}
  \label{sz2aw} x \overline{r}_n(x)=\overline{r}_{n+1}(x)+(\alpha_A A_n+ \alpha_C
C_n)\overline{r}_n(x)+\alpha_A\alpha_C A_{n-1}C_n
\overline{r}_{n-1}(x),
\end{equation}
To see that  \eqref{sz2aw}
is equivalent to the Askey-Wilson recurrence \eqref{AW1} we only need to show that $\alpha_A=\alpha_C$.
Equations \eqref{EqC} give
$$a=\eta\frac{\sqrt{(1-q)t+\tau}}{\sqrt{(1-q)^2+\eta\theta(1-q)+\eta^2\tau}}.$$
Inserting this into the expressions for $\alpha_A$ and $\alpha_C$,
we see that
$$
\alpha_A=\alpha_C=\frac{{\sqrt{(1-q)t+\tau}}{\sqrt{(1-q)^2+\eta\theta(1-q)+\eta^2\tau}}}{(1-q)^2}.
$$

 Our main result is the
following projection formula.
\begin{theorem}\label{T1} If $0\leq s\leq t$, $0\leq q\leq 1$, $\eta\theta\geq 0$, $\tau\geq0$,
 then
for all $x$ in the support of 
the orthogonality measure
$\mu_s$ there exists a unique
probability measure $\nu_x=\nu_{x,t,s}$ such that
\begin{equation}\label{proj main}
\overline{p}_n(x;s)=\int \overline{p}_n(y;t)\nu_x(dy).
\end{equation}

\end{theorem}
Of course, probability measure $\nu_x=\nu_{x,t,s}$ depends also on
parameters $0\leq s\leq t$ as well as on the remaining parameters
$\eta,\theta,\tau,q$.

\begin{remark}
Since Askey-Wilson polynomials are symmetric with respect parameters $b,c,d$,
\cite[Formula (2.4)]{Rahman85} can be written as
$$\int \bar{w}_n(y;a,b,c,d)\nu_x(dy)=\bar{w}_n(x;\mu a, \mu^{-1}b,\mu^{-1}c,\mu d).$$
Under this parametrization, from \eqref{EqA} and \eqref{EqC} e can see that
Theorem \ref{T1} is essentially this formula with $d=0$,
 and $\mu^2=((1-q)s+\tau)/(1-q)t+\tau) \leq 1$. The only improvement over
 \cite{Rahman85} is that
 due to our interest in applications to probability, our range of parameters
 leads to measures $\nu_x(dy)$ that might have a discrete component.
\end{remark}

The proof of Theorem \ref{T1} appears in Section \ref{Proof of T1},
after a number of preliminary results. The plan of the proof is as
follows.   In Section \ref{Aux Poly} we define a family of monic
polynomials $\{\overline{Q}_n\}$ in variable $y$. We verify that the
assumptions of Favard's theorem are satisfied for the relevant pairs
$(x, s)$, so that their orthogonality  measure $\nu_{x,t,s}$ exists.
We show that  this measure is unique (a fact that is nontrivial only
when $q=1$).  We then use the formula for the connection
coefficients between polynomials $\{\overline{Q}_n\}$ and
 $\{\overline{p}_n\}$ to deduce \eqref{proj main}.

When $\eta>0$, we will find it convenient to consider the following
non-monic polynomials
\begin{equation}
\label{p-rec} yp_n(y;t)=a_n p_{n+1}(y;t)+b_n(t) p_n(y;t)+c_n(t)
p_{n-1}(y;t).
\end{equation}
Clearly they have the same orthogonality measure $\mu_t$ as the
monic polynomials.

\section{Identities}

We will need a number of auxiliary identities.
\begin{lemma}\label{L1}
Fix a sequence $\{p_n: n\geq 0\}$ of real  numbers. Let $\{\beta_{n,k}:0\leq k \leq
n, n=0,1,\dots\}$ be defined   by  $\beta_{n,k}=0$ for $n<0$ or
$k>n$ and for $0\leq k\leq n$ by the recurrence
\begin{equation}\label{beta def}
 [k+1]_q\beta_{n,k+1}=q^k[n-k]_q \beta_{n,k}+[n]_q\beta_{n-1,k},
\end{equation}
with the initial values $\beta_{n,0}=p_n$, $n=0,1,\dots$. Then
\begin{equation}\label{beta sol}
\beta_{n,k}=\left[ \begin{array}{c}n \\ k \end{array}\right]_{q}
\sum_{j=0}^k\left[ \begin{array}{c}k \\ j \end{array}\right]_{q}
q^{(k-j)(k-j-1)/2}p_{n-j}.
\end{equation}
\end{lemma}
\begin{proof} This follows by a routine induction argument with respect to $k$.
Clearly \eqref{beta sol} holds true for $k=0$ and all $n\geq 0$.
Suppose \eqref{beta sol} holds true for some $k\geq 0$ and all
$n\geq 0$. Then by \eqref{beta def} and the induction assumption, we
have
\begin{multline*}
\beta_{n,k+1}=\frac{[n-k]_q}{[k+1]_q}\left[ \begin{array}{c}n \\ k
\end{array}\right]_{q}\sum_{j=0}^k
\left[ \begin{array}{c}k \\ j \end{array}\right]_{q}
q^{k+(k-j)(k-j-1)/2} p_{n-j}+\\\frac{[n]_q}{[k+1]_q}\left[
\begin{array}{c}n-1 \\ k
\end{array}\right]_{q}\sum_{j=0}^k\left[ \begin{array}{c}k \\ j \end{array}\right]_{q}
q^{(k-j)(k-j-1)/2}  p_{n-1-j}\\ =\left[ \begin{array}{c}n \\ k+1
\end{array}\right]_{q} p_{n-(k+1)}+ \left[ \begin{array}{c}n \\ k+1
\end{array}\right]_{q} q^{(k+1)k/2} p_{n}\\+\left[
\begin{array}{c}n \\ k+1 \end{array}\right]
_{q}\sum_{j=1}^k\left(q^j\left[ \begin{array}{c}k \\ j
\end{array}\right]_{q}+\left[ \begin{array}{c}k \\ j-1 \end{array}\right]_{q}
\right)q^{(k+1-j)(k-j)/2}p_{n-j}.
\end{multline*}
The well known formula \cite[(I.45)]{Gasper-Rahman-90}
\begin{equation}\label{WellKnown}
\left[ \begin{array}{c}k+1 \\ j \end{array}\right]_{q} =q^j\left[
\begin{array}{c}k \\ j \end{array}\right]_{q}+\left[
\begin{array}{c}k \\ j-1 \end{array}\right]_{q}=\left[
\begin{array}{c}k \\ j \end{array}\right]_{q}+q^{k-j+1}\left[
\begin{array}{c}k \\ j-1 \end{array}\right]_{q},\; 1\leq j\leq k,
\end{equation}
ends  the proof.
\end{proof}
It turns out that expressions of the form \eqref{beta sol} can
sometimes be written as products.
\begin{proposition}\label{PP3}
If polynomials $\{p_n(y;t)\}$ satisfy recurrence \eqref{p-rec} and
\begin{eqnarray}\label{AA}
  A_n(x,s)&=& a_n+q^n x - s\eta q^n[n]_q\;,
\end{eqnarray}
then
for all $k\geq 1$ we have
\begin{equation}\label{**}
\sum_{j=0}^k\left[ \begin{array}{c}k \\ j \end{array}\right]_{q}
q^{(k-j)(k-j-1)/2}p_{k-j}(x;s)=
\prod_{j=0}^{k-1}\frac{A_j(x,s)}{a_j}.
\end{equation}
\end{proposition}
\begin{proof} We proceed by induction with respect to $k$. Formula
\eqref{**} holds true for $k=0$ by convention, and for $k=1$ by a
calculation: $p_1(x;s)+p_0(x;s)=1+\eta x$.

Let $\beta_{n,k}(x,s)$ be defined by \eqref{beta sol} with
$p_n=p_n(x;s)$, $n=0,1,\dots$. The induction assumption says that
\begin{equation}
  \label{*1}
  \beta_{k,k}(x,s)=\prod_{j=0}^{k-1}\frac{A_j(x,s)}{a_j}
\end{equation} for some
$k\geq 1$. From  \eqref{beta def} we see that
\begin{equation}\label{*2}
  \beta_{k+1,k+1}(x,s)=\beta_{k,k}(x,s)+q^k\sum_{j=0}^k\left[ \begin{array}{c}k \\ j \end{array}\right]_{q}
q^{(k-j)(k-j-1)/2}p_{k+1-j}(x;s).
\end{equation}
On the other hand, multiplying both sides of \eqref{*1} by
$\frac{A_k(x,s)}{a_k}=1+\frac{x-s\eta[k]_q}{a_k}q^k$ and using
\eqref{p-rec} we see that
\begin{multline}\label{*3}
\prod_{j=0}^{k}\frac{A_j(x,s)}{a_j}=\beta_{k,k}(x,s) -q^k\frac{s\eta
[k]_q}{a_k}\sum_{j=0}^k \left[
\begin{array}{c}k \\ j \end{array}\right]_{q}
q^{(k-j)(k-j-1)/2}p_{k-j}(x;s)\\+\frac{q^k}{a_k}\sum_{j=0}^k \left[
\begin{array}{c}k \\ j \end{array}\right]_{q}
q^{(k-j)(k-j-1)/2}\big(a_{k-j}p_{k+1-j}(x;s)+b_{k-j}(s)p_{k-j}(x;s)\\+c_{k-j}(s)p_{k-1-j}(x;s)\big).
\end{multline}
Writing the right hand side of \eqref{*3} as
\begin{equation}
  \beta_{k,k}(x,s) +\frac{q^k}{a_k}\sum_{j=0}^{k+1} \left[
\begin{array}{c}k \\ j \end{array}\right]_{q}
q^{(k-j)(k-j-1)/2}\gamma_{k,j}p_{k+1-j}(x;s),
\end{equation}
from \eqref{aa}, \eqref{bb} and \eqref{cc} it is not difficult to
see that $\gamma_{k,0}=a_k$ and $\gamma_{k,k+1}=-s\eta[k]_q+s\eta
[k]_q =0$. Similarly, for $1\leq j \leq k$ we have
\begin{multline}
  \label{*4}
  \gamma_{k,j}=a_{k-j}+ \frac{q^{k-j}[j]_q}{[k+1-j]_q}\left(b_{k+1-j}(s)-s\eta[k]_q\right)\\+
  \frac{q^{2k-2j+1}[j]_q[j-1]_q}{[k+2 -j]_q[k+1-j]_q}c_{k+2-j}(s)
  =\eta^{-1}+\theta[k-j]_q+\eta\tau[k-j]_q^2+[j]_q\theta
  q^{k-j}\\+\eta\tau \Big([j]_q[k+1-j]_qq^{k-j}+[j]_q[k-j]_qq^{k-j}+[j]_q[j-1]_qq^{2k+1-2j}\Big)
 \\ =
  \eta^{-1}+\theta\left([k-j]_q+q^{k-j}[j]_q\right)\\+\eta\tau\left([k-j]_q\left([k-j]_q+q^{k-j}[j]_q\right)+  [j]_q
  q^{k-j}\left([k+1-j]_q+[j-1]_qq^{k+1-j}\right)\right)\\
  =\eta^{-1}+\theta[k]_q+\eta\tau[k]_q^2.
\end{multline}
(Here we used repeatedly the identity $[k-j]_q+q^{k-j}[j]_q=[k]_q$.) Thus $\gamma_{k,j}=a_k$,
which shows that the right hand sides of equations \eqref{*2} and
\eqref{*3} are equal. Therefore their left hand sides are equal,
ending the proof.
\end{proof}
 For $n\geq 0$, $\eta\ne 0$,  and $q\ne 0$ let
\begin{equation}\label{x}
x_n(s)=\frac{(s \eta
q^n-\theta)[n]_q-\eta^{-1}-\eta\tau[n]_q^2}{q^n}\end{equation} be
the zero of $A_n(x,s)$, see \eqref{AA}. It turns out that \eqref{**}
extends to higher order polynomials $p_n$ when the polynomials are
evaluated at $x_k$.

\begin{lemma}\label{L2}
If $\{p_n(y;t)\}$ satisfies recurrence \eqref{p-rec} and $\eta, q>
0$ then for $n\geq k$ we have
\begin{equation}\label{*}
\sum_{j=0}^k\left[ \begin{array}{c}k \\ j \end{array}\right]_{q}
q^{(k-j)(k-j-1)/2}p_{n-j}(x_k(s);s)=
\frac{(-1)^{n-k}}{q^{k(n-k)}}\prod_{j=0}^{k-1}\frac{A_j(x_k(s),s)}{a_j}.
\end{equation}
(For $k=0$ this should be interpreted as $p_n(x_0;s)=(-1)^n$, $n\geq
0$.)
\end{lemma}

\begin{proof}  Let $\beta_{n,k}(x,s)$ be defined by \eqref{beta sol} with
$p_n=p_n(x;s)$, $n=0,1,\dots$. Then the left hand side of \eqref{*}
is $\beta_{n,k}(x_k(s),s)/\left[ \begin{array}{c}n \\ k
\end{array}\right]_{q}$. We first prove an auxiliary fact that for
all $0\leq j<k\leq n$ we have $\beta_{n,k}(x_j(s),s)=0$. We prove
this by induction with respect to $n-k$. Suppose there is $m\geq 0$
such that $\beta_{n,k}(x_j(s),s)=0$ for all triplets $(j,k,n)$ such
that $0\leq j<k$ and $n-k= m$. By \eqref{**} this holds true for
$m=0$. Given $j<k$ and $n$ such that $n-k=m+1$ by \eqref{beta def}
we have
\begin{equation}
([n]_q-[k]_q)\beta_{n,k}(x,s)=[k+1]_q
\beta_{n,k+1}(x,s)-[n]_q\beta_{n-1,k}(x,s). \label{!}
\end{equation}
 By
induction assumption the right hand side of \eqref{!} evaluated at
$(x_j(s),s)$ vanishes. As $q\ne 0$ and $n=k+m+1>k$, we have
$[n]_q-[k]_q\ne 0$, so $\beta_{n,k}(x_j(s),s)=0$.

We now prove \eqref{*}. From \eqref{sum=1/eta} it is easy to see by
induction that $p_n(x_0;s)=p_n(-\eta^{-1};s)=(-1)^n$. For $k\geq 1$,
we will prove \eqref{*} by induction with respect to $n$.

 If $n=k$ then formula \eqref{*} holds by Proposition \ref{PP3}.
Suppose \eqref{*} holds for some $n\geq k$. Then from \eqref{beta
def} and the fact that $\beta_{n+1,k+1}(x_k(s),s)=0$ we see that
\begin{multline*}
\beta_{n+1,k}(x_k(s),s)=-\frac{[n+1]_q}{q^k[n+1-k]_q}\beta_{n,k}(x_k(s),s)\\=-\frac{[n+1]_q}{q^k[n+1-k]_q}\left[
\begin{array}{c}n \\ k
\end{array}\right]_{q} \frac{(-1)^{n-k}}{q^{k(n-k)}}\prod_{j=0}^{k-1}\frac{A_j(x_k(s),s)}{a_j}.
\end{multline*}
Therefore, \
\begin{multline*}
\sum_{j=0}^k\left[ \begin{array}{c}k \\ j \end{array}\right]_{q}
q^{(k-j)(k-j-1)/2}p_{n+1-j}(x_k(s);s)=\frac{\beta_{n+1,k}(x_k(s),s)}{\left[
\begin{array}{c}n+1 \\ k
\end{array}\right]_{q}}\\=\frac{(-1)^{n+1-k}}{q^{k(n+1-k)}}\prod_{j=0}^{k-1}\frac{A_j(x_k(s),s)}{a_j}.
\end{multline*}
\end{proof}
We need to analyze equation \eqref{*} in more detail.

\begin{lemma}\label{l2.5} Fix $k\geq 1$, $q\neq 0$, $\Pi_k>0$. Suppose that $(p_n)_{n\geq 0}$ is the general solution of the recurrence
\begin{equation}\label{inversion}
\sum_{j=0}^k\left[ \begin{array}{c}k \\ j \end{array}\right]_{q}
q^{(k-j)(k-j-1)/2}p_{n-j}= \frac{(-1)^{n-k}}{q^{k(n-k)}}\Pi_k,\;
n\geq k.
\end{equation} Then
\begin{equation}\label{inversion sol}
p_n=(-1)^{n-k}q^{-nk}q^{k(k+1)/2}  \left(\left[\begin{array}{c}n
\\ k \end{array}\right]_{q}+\sum_{r=1}^k
C_r q^{nr}\left[\begin{array}{c}n \\ k-r \end{array}\right]_{q}\right)\Pi_k, \; n\geq 0,
\end{equation}
where $C_1,\dots,C_k$ are arbitrary constants.

\end{lemma}
\begin{proof}
Substitute $$y_n=\frac{(-1)^{n-k} q^{k(2n- k-1)/2}}{\Pi_k } p_n.$$
Then with $y=(y_n)_{n\geq 0}$ the equation takes the form of an
initial value problem for a linear recurrence with constant
coefficients:
\begin{eqnarray}
  \label{Delta_q}
  (\Delta_{q,k} y)_n&=& 1, \; n\geq k,
\end{eqnarray}
where
$$(\Delta_{q,k} y)_n= \sum_{j=0}^k(-1)^j\left[ \begin{array}{c}k \\ j \end{array}\right]_{q}
q^{j(j+1)/2} y_{n-j}.$$

We remark that when $q=1$ we trivially have
$\Delta_{1,k}=\Delta_{1,1}^k$. Since $\Delta_{1,1} $ is the usual
difference operator, in this case the general solution of
\eqref{Delta_q} is well known. The  $q$-generalization of
this formula follows from \eqref{WellKnown} by induction with
respect to $k$. We have
\begin{equation}\label{Factor Delta}
\Delta_{q,k}=R_1R_2\dots R_k,
\end{equation}
where $R_j=\Delta_{q^j,1}$ are commuting difference operators, $
(R_j y)_n=y_n-q^jy_{n-1}$ for $n\geq 1$.

The general theory of linear difference equations implies that
\eqref{inversion sol} is a consequence of the following two observations.
\begin{claim}\label{claim sol NH}
\begin{enumerate}
\item $(\Delta_{q,k} y)_n=1$ for $n\geq k$ when
\begin{equation}\label{x=1} y_n= \left[
\begin{array}{c}n  \\ k \end{array}\right]_{q},
\end{equation}
\item  $(\Delta_{q,k} y)_n=0$ for $n\geq k-r$ when \begin{equation}
\label{x=0}y_n= \left[
\begin{array}{c}n \\ k-r \end{array}\right]_{q}  q^{r n} , \; r=1,2,\dots,k.
\end{equation}
\end{enumerate}
\end{claim}
\begin{proof}[Proof of Claim \ref{claim sol NH}]
We note that \eqref{x=1} is just $r=0$ case of \eqref{x=0}.

For fixed $r\geq 0$ and $n\geq k\geq r$ we have
\begin{multline*}
\left(R_k \left(\left[
\begin{array}{c}n \\ k-r \end{array}\right]_{q}  q^{r n} \right)\right)_n=
\left[
\begin{array}{c}n \\ k-r \end{array}\right]_{q}  q^{r n}-\left[
\begin{array}{c}n-1 \\ k-r \end{array}\right]_{q}  q^{r n+k-r}\\
=q^{rn}\frac{[n-k+r+1]_q\dots
[n-1]_q}{[k-r]_q!}\left([n]_q-q^{k-r}[n-k+r]_q\right)\\
=q^{rn}\frac{[n-k+r+1]_q\dots [n-1]_q}{[k-r]_q!}[k-r]_q\\
=\left[
\begin{array}{c}n -1\\ k-1-r \end{array}\right]_{q}  q^{r n}.
\end{multline*}
Therefore
$$
\left(R_{r+1}R_{r+2}\dots R_k \left(\left[
\begin{array}{c}n \\ k-r \end{array}\right]_{q}  q^{r n}
\right)\right)_n=\left[
\begin{array}{c}n-k+r \\ 0 \end{array}\right]_{q}  q^{r n}=q^{rn}.
$$
If $r=0$ this implies \eqref{x=1} by \eqref{Factor Delta}. If $r\geq
1$ then to prove \eqref{x=0} it remains to notice that since $n\geq
k\geq r \geq 1$ we have $\left(R_r(q^{nr})\right)_n=q^{rn}-q^r q^{(n-1)r}=0$.
\end{proof}
The constants $C_1,\dots,C_k$ are determined from the condition that
formula \eqref{inversion sol} holds for $p_0,\dots,p_{k-1}$.
\end{proof}

\begin{proposition}\label{PP4.5}
Suppose $\{p_n(y;t)\}$ satisfies recurrence \eqref{p-rec}. Then
there are constants $c_k(s)$ that do not depend on $n$ such that:
\begin{enumerate}
\item if $0<q<1$ then
$|p_n(x_k(s);s)|\leq c_k(s)q^{-kn}$;
\item if $q=1$ then $|p_n(x_k(s);s)|\leq c_k(s) n^{k}$.
\end{enumerate}
\end{proposition}
\begin{proof}
This follows from \eqref{inversion sol} and \eqref{*}.
\end{proof}

\section{Uniqueness of the moment problem}
\begin{proposition}\label{TT} Suppose $0\leq q\leq 1$, $\eta>0$, $\theta\geq 0$,
$\tau\geq 0$.  Let  $\{p_n(y;t)\}$  be defined by \eqref{p-rec}.
Then the orthogonality measure $\mu_t$ of polynomials $\{p_n(y;t)\}$
is determined uniquely by moments.
\end{proposition}

\begin{proof} For $|q|<1$, the coefficients of the recurrence are
bounded, so the only case that requires proof is $q=1$. Furthermore,
the conclusion holds for $\tau=0$, as in this case $\mu_t$ is a
negative binomial law, see \cite{Bryc-Wesolowski-05}. It therefore
remains to consider the case   $q=1$, $\tau>0$.

 In this case, we use
the fact that with $x_0=-\eta^{-1}$ we have $p_n(x_0;t)=(-1)^n$, see
Lemma \ref{L2}. Let $q_n(y;t)$ be the associated polynomials which
satisfy recurrence \eqref{p-rec} for $n\geq 1$ with the initial
terms $q_0=0$, $q_1=1/a_0$. Then
$$x_0q_n(x_0)=a_nq_{n+1}(x_0)+(a_n+c_n +x_0)q_{n}(x_0)+c_nq_{n-1}(x_0).$$
Therefore with $f_n(t):=(-1)^{n-1}q_n(x_0;t)$ we have
\begin{multline*}
f_{n+1}-f_n=\frac{c_n(t)}{a_n}\left(f_{n}-f_{n-1}\right)=\frac{c_1(t)c_2(t)\dots
c_n(t)}{a_1a_2\dots a_n}\left(f_1-f_0\right)\\=
\frac{c_1(t)c_2(t)\dots c_n(t)}{a_0a_2\dots a_{n-1}}\frac{1}{a_n}.
\end{multline*}
Thus 
with a suitable convention for $n=1$ we can write the solution as
\begin{equation}
  \label{f}
  f_{n+1}(t)=\sum_{k=0}^n\frac{c_1(t)c_2(t)\dots c_k(t)}{a_0a_2\dots
a_{k-1}}\;\frac{1}{a_k}.
\end{equation}
 Let
\begin{equation}\label{p orthonormal}
\widetilde{p}_n(x;t)=\sqrt{\frac{a_0a_2\dots
a_{n-1}}{c_1(t)c_2(t)\dots c_n(t)}}\;p_n(x;t) \end{equation} and
$$\widetilde{q}_n(x;t)=\sqrt{\frac{a_0a_2\dots
a_{n-1}}{c_1(t)c_2(t)\dots c_n(t)}}\;q_n(x;t)
$$ be the corresponding orthonormal polynomials.

By \cite[page 84]{Akhiezer65}, the moment problem is determined
uniquely, if
\begin{equation}\label{infinity}
\sum_n|\widetilde{p}_n(x_0)|^2+\sum_n|\widetilde{q}_n(x_0)|^2=\infty.
\end{equation}
We have
$$
|\widetilde{p}_n(x_0)|^2=\frac{a_0a_2\dots
a_{n-1}}{c_1(t)c_2(t)\dots c_n(t)},
$$
and from \eqref{f} we get
$$
|\widetilde{q}_n(x_0)|^2=\frac{a_0a_2\dots
a_{n-1}}{c_1(t)c_2(t)\dots
c_n(t)}\left(\sum_{k=0}^n\frac{c_1(t)c_2(t)\dots c_k(t)}{a_0a_2\dots
a_{k-1}}\;\frac{1}{a_k}\right)^2.
$$
To verify \eqref{infinity} we use the fact that
\begin{equation}
  \label{elementary estimates}
a_n\approx \eta\tau n^2,\;
\frac{c_{n+1}(t)}{a_n}=1+\frac{\alpha(t)}{n}+O(1/n^2),\;\frac{a_n}{c_{n+1}(t)}=1
-\frac{\alpha(t)}{n}+O(1/n^2),\end{equation} where
\begin{equation}
  \label{alpha(t)}
  \alpha(t)= \frac{t\eta-\theta}{\eta\tau}+1,
\end{equation}
and $a_n\approx b_n$ means that $a_n/b_n\to 1$.

If $t\leq \theta/\eta$ then $\alpha(t)\leq 1$ and
$$\prod_{k=1}^n\left(1-\frac{\alpha(t)}{k}\right)\approx
\exp\left(-\alpha(t)\sum_{k=1}^n\frac1k\right)\approx
n^{-\alpha(t)}\geq n^{-1},$$ so the first series in \eqref{infinity}
diverges. On the other hand, if $t>\theta/\eta$ so that
$\alpha(t)>1$, then
$$
|\widetilde{q}_n(x_0)|^2\approx
n^{-\alpha(t)}\left(\sum_{k=1}^{n-1}k^{\alpha(t)}\frac{1}{k^2}\right)^2\approx
n^{-\alpha(t)}\left(n^{\alpha(t)-1}\right)^2=
n^{\alpha(t)-2}>n^{-1},
$$
so the second series in \eqref{infinity} diverges.
\end{proof}

\section{Support of the orthogonality measure}
 Recall
that $\mu_t$ denotes the orthogonality measure of polynomials
$\{p_n(y;t)\}$. The following result will be used to define the
orthogonality measure of auxiliary polynomials in Section \ref{Aux
Poly}.
\begin{proposition}\label{P+} Suppose $0\leq q\leq 1$, $\eta>0$, $\theta\geq 0$,
$\tau\geq 0$. If $x\in\supp(\mu_s)$, then
\begin{equation}
  \label{Prod A}
  \prod_{j=0}^nA_j(x,s)\geq 0\mbox{ for all }n\geq 0.
\end{equation}
\end{proposition}

We prove Proposition \ref{P+} from rudimentary information about the
support of $\mu_s$. 

\begin{lemma}\label{T supp}
Let
$x_j(t)$ be given by \eqref{x}. Then the support of $\mu_t$ is a
subset of the interval $[x_0(t),\infty)$. In addition, if
$t>\frac{\theta}{\eta}+\frac{1-q}{\eta^2}$ then
$$\supp(\mu_t)\subset \{x_0(t),
x_1(t),\dots,x_{k_*}(t)\}\cup[y_*,\infty)$$ with
$y_*=\max\{x_{k_*}(t), x_{k_*+1}(t)\}$ and
\begin{equation}\label{k=}
k_*=\begin{cases} \max\left\{k:t>
\frac{\theta}{\eta}+2\tau k\right\},& q=1;\\
\\ \displaystyle
\max\left\{k: q^{2k}>
\frac{(1-q)^2+(1-q)\eta\theta+\eta^2\tau}{\eta^2(t(1-q)+\tau)}
\right\},& 0<q<1.
\end{cases}
\end{equation}
\end{lemma}
\begin{remark} We note that $y_*=x_{m^*}(t)$ with
\begin{equation*}
m^*=\begin{cases} \max\{j:t>
\frac{\theta}{\eta}+2\tau j-1\}, & q=1;\\
\max\{j: q^{2j-1}\eta^2(t(1-q)+\tau)>
(1-q)^2+(1-q)\eta\theta+\eta^2\tau \}, & q<1.
\end{cases}
\end{equation*}
\end{remark}
We will use the following criterion to show that there are at most
$k_*+1$ atoms below $y_*$.

\begin{CT}\label{Sturm} Suppose $p_n(x)$ are orthogonal polynomials with unique
orthogonality measure $\mu$. If the sequence $\{(-1)^np_n(a): n\geq
0\}$ changes sign $k$-times, then there is a finite set  $D$ with at most $k$ points such that
 $$\supp(\mu_t)\subset D\cup[a,\infty).$$
In particular, $\mu$ has at most $k$ atoms in $(-\infty,a)$.
\end{CT}
\begin{proof}This follows  from the
interlacing property of zeros of orthogonal polynomials. The details
are omitted.
\end{proof}
\begin{proof}[Proof of Lemma \ref{T supp}]
We first observe that $\supp(\mu_s)\subset[-\eta^{-1},\infty)$. This
follows from the fact that by Proposition \ref{TT} measure $\mu_s$
is determined uniquely, so we can combine Lemma \ref{L2} applied to
$k=0$ with Theorem \ref{Sturm} applied to $a=x_0=-\eta^{-1}$.

 We now verify that if $t>\frac{\theta}{\eta}+\frac{1-q}{\eta^2}$
then there are $k_*+1$ atoms at $\{x_0(t),
x_1(t),\dots,x_{k_*}(t)\}$. Recall that $x_j(t)$ is an atom of
$\mu_t$ if the orthonormal polynomials \eqref{p orthonormal} are
square-summable at $x=x_j(t)$, see \cite[page 84]{Akhiezer65}. We
will consider separately the cases $q=1$ and $0<q<1$.

Suppose $q=1$. Then by Proposition \ref{PP4.5} we have
$$\sum_n|\widetilde{p}_n(x_j;t)|^2\leq c_j \sum_n n^{2j}\frac{a_0a_2\dots
a_{n-1}}{c_1(t)c_2(t)\dots c_n(t)}\approx c_j\sum_n n^{2j-\alpha(t)}. $$
(See \eqref{alpha(t)}.) Therefore from \eqref{elementary estimates},
the series converges if $t>\frac{\theta}{\eta}+2j\tau$.

Suppose now that $0<q<1$. Then by Proposition \ref{PP4.5} we have
$$\sum_n|\widetilde{p}_n(x_j;t)|^2\leq c_j \sum_n q^{-2nj}\frac{a_0a_2\dots
a_{n-1}}{c_1(t)c_2(t)\dots c_n(t)}. $$ Since
$$\lim_{n\to\infty} \frac{a_n}{c_{n+1}(t)}=\frac{(1-q)^2+(1-q)\eta\theta +\eta^2\tau}{\eta^2(t(1-q)+\tau)},$$
the series converges if $\frac{(1-q)^2+(1-q)\eta\theta
+\eta^2\tau}{q^{2j}\eta^2(t(1-q)+\tau)}<1$. This proves that
$x_j(t)$, $0\leq j\leq k_*$ is an atom under the condition
\eqref{k=}.

To estimate that there are at most $k_*+1$ atoms below $y_*$ we use
Lemma \ref{l2.5} to verify that there are at most $k$ atoms of
$\mu_t$ below $x_k(t)$. Namely, Lemma \ref{l2.5} states that there
exists a polynomial $r(x)$ of degree $k$ such that
$$p_n(x_k(t);t)=\begin{cases}
r(n), & q=1;\\
r(q^{-n}),&0<q<1.
\end{cases}
$$
Since $r(x)=0$ has at most $k$ real solutions, the sequence
$\{p_n(x_k(t);t):n\geq 0\}$ has at most $k$ changes of sign.
Proposition  \ref{TT} implies that we can use Theorem \ref{Sturm} to
end the proof.
\end{proof}

\begin{proof}[Proof of Proposition \ref{P+}]
 If $q=0$ then $A_n(x,s)=\eta^{-1}$ does not depend on $x$ for $n\geq 1$.
Since $A_0(x,s)=1/\eta+\theta+\eta\tau+x$, \eqref{Prod A} follows
from $\supp(\mu_s)\subset[-\eta^{-1},\infty)$.

In the remaining part of the proof, we assume $0<q\leq 1$. We use
the trivial observation that $A_j(x,s)$ increases as a function of
$x$ and decreases as a function of $s$.

 Suppose $0\leq s\leq \theta/\eta+(1-q)/\eta^2$. From
$x\in \supp(\mu_s)\subset[-\eta^{-1},\infty)$ we get
$$A_n(x,s)\geq A_n\left(-\frac{1}{\eta},\frac{\theta}{\eta}+\frac{1-q}{\eta^2}\right)=
\frac{(1-q^n)^2}{\eta}+\theta [n]_q(1-q^n)+[n]_q^2\eta\tau\geq 0. $$
Thus \eqref{Prod A} holds.

Suppose $s> \theta/\eta+(1-q)/\eta^2$ so that $k_*=k_*(s)\geq 0$ is
well defined. We notice that
\begin{equation}\label{x-ineq}
x_0(s)\leq x_1(s)\leq \dots \leq  x_{k_*}(s)\leq y_* \mbox{ and }
x_j(s)\leq  y_* \mbox{ for all } j>k_*.
\end{equation}
Omitting the easier case of $q=1$, write $x_n(s)=h(q^n)$, where
$$
h(z)=-\frac{1}{\eta z } + \frac{\left( 1 - z
\right)  \left( s z \eta  - \theta  \right) }{z(1 - q)} -
    \frac{{\left( 1 - z \right) }^2 \eta  \tau }{{z\left( 1 - q \right)
    }^2}.
$$
A calculation shows that
$$
h''(z)=-2\frac{ \left( 1 - q \right)  \left( 1 - q + \eta  \theta
\right)  + {\eta }^2 \tau  }
  {{\left( 1 - q \right) }^2 z^3 \eta }<0
$$  on the interval $0<z<1$.
Since $h$ tends to $-\infty$ at the endpoints, therefore it has a
unique maximum $z_*\in(0,1)$ given by
$$z_*^2
    {\eta }^2 \left( \left( 1 - q \right)  s + \tau  \right) = \left( 1 - q \right)  \left( 1 - q + \eta
 \theta  \right)  + {\eta }^2 \tau
  .$$
 In particular, $q^{k_*+1} \leq z_*< q^{k_*}$, so $h(z)$ increases on $(0,q^{k_*+1})$ and decreases on
$(q^{k_*},1]$. Thus
 $h(q^{k_*+1})\geq
h(q^{k_*+2})\geq\dots$, and $h(q^0)\leq h(q^1)\leq \dots \leq
h(q^{k_*})$.

Inequality \eqref{x-ineq} ends the proof as follows. If $x=x_k(s)$
for some $k\leq k_*$ then $A_j(x,s)\geq A_j(x_j(s),s)=0$ for $0\leq
j\leq k$, so \eqref{Prod A} holds for $0\leq n<k$. On the other
hand, $A_k(x,s)=0$, so \eqref{Prod A} holds trivially for all $n\geq
k$.

Suppose now that $x\geq y_*$. Then \eqref{x-ineq} implies
$A_j(x,s)\geq A_j(y_*,s)\geq A_j(x_j(s),s)=0$ for all
$j=0,1,2\dots$. Thus \eqref{Prod A} follows.

\end{proof}

\section{Auxiliary polynomials}\label{Aux Poly} For the proof of
Theorem \ref{T1} we construct measure $\nu$ as a measure of
orthogonality of auxiliary monic polynomials $\overline{Q}_n(y;x,t,s)$ in
variable $y$. We begin with a non-monic version of these
polynomials,  defined by the three step recurrence
\begin{multline}\label{Q rec}
  y\  Q_n(y;x,t,s)=A_n(x,s)Q_{n+1}(y;x,t,s)+B_n(x,t,s) Q_n(y;x,t,s)\\+C_n(t,s)
  Q_{n-1}(y;x,t,s),
\end{multline}
where $A_n$ is defined by \eqref{AA} and
\begin{eqnarray}
 B_n(x,t,s)&=& b_n(t)+q^n x-(1+q)q^{n-1}\eta s[n]_q\;,
 \label{BB} \\
 C_n(t,s)&=& c_n(t)-q^{n-1}\eta s[n]_q\;.
 \label{CC}
\end{eqnarray}
with $Q_{-1}=0$, $Q_0=1$. The Jacobi matrix of this recurrence
arises as a solution of the $q$-commutation equation
\cite[(1)]{Bryc-Matysiak-Wesolowski-04} with the appropriately
modified initial condition; for more details see
\cite{Bryc-Matysiak-Wesolowski-05}.
Polynomials $\{Q_n\}$ are well defined for all $x,s,t$ as long as
$x\not\in\{x_0(s),x_1(s),\dots\}$.

\subsection{Connection Coefficients}
For $x\not\in\{x_0(s),x_1(s),\dots\}$, the connection coefficients
$\beta_{n,k}(x,t,s)$ are defined implicitly by
\begin{equation}\label{p2Q}
  p_n(y;t)=\sum_{k=0}^n \beta_{n,k}(x,t,s)Q_k(y;x,t,s).
\end{equation}
Our next goal is to find the connection coefficients
$\beta_{n,k}(x,t,s)$ explicitly and to show that they do not depend
on $t$.

Define two linear operators $K,L:\sR^\infty\to\sR^\infty$  acting on
infinite matrices  $\beta=[\beta_{n,k}]_{n,k\geq 0}$ by the rule
$$
[K\beta]_{n,k}=a_n \beta_{n+1,k}+b_n(t)
\beta_{n,k}+c_n(t)\beta_{n-1,k},
$$
$$
[L\beta]_{n,k}=A_{k-1}(x,s) \beta_{n,k-1}+B_k(x,t,s)
\beta_{n,k}+C_{k+1}(t,s)\beta_{n,k+1}.
$$
Let
\begin{equation}\label{HH}
H_t\beta=K\beta-L\beta. \end{equation}
 The infinite triangular
matrix $[\beta_{n,k}(x,t,s)]_{n\geq k\geq 0}$ is a unique solution
of the discrete boundary value problem
\begin{equation}
  \label{H}
  H_t\beta=0,
\end{equation}
\begin{equation}\label{H ini}
\beta_{n,n}(x,t,s)=\prod_{j=0}^{n-1}\frac{A_j(x,s)}{a_j}, \;n\geq 0.
\end{equation}
The boundary condition \eqref{H ini} arises from \eqref{p2Q} by
comparing the coefficients at $y^n$. Equation \eqref{H} follows
directly from the recurrences; here we give an argument based on the
fact that the orthogonality measure for polynomials $\{Q_n\}$ exists
for an infinite set of $x$. For such $x$, we have
$$\beta_{n,k}(x,t,s)=\frac{\int p_n(y;t)Q_k(y;x,t,s)\nu_{x,t,s}(dy)}{\|Q_k\|_2^2}.$$
Since
$$\|Q_k\|_2^2=\prod_{j=1}^{k}\frac{C_j(t,s)}{A_{j-1}(x,s)},$$
 \eqref{H} follows from
$$\int [yp_n(y;t)]Q_k(y;x,t,s)\nu_{x,t,s}(dy)=\int
p_n(y;t)[yQ_k(y;x,t,s)]\nu_{x,t,s}(dy).
$$ by \eqref{p-rec} and
\eqref{Q rec}. Of course, once \eqref{H} holds for a large enough
set of $x$, it holds for all $x$.
\begin{lemma}\label{L5}
If $x\not\in\{x_0(s),x_1(s),\dots\}$ then the coefficients
$\beta_{n,k}(x,t,s)$ in \eqref{p2Q} are determined uniquely, and do
not depend on variable $t$. In fact,
$\beta_{n,k}(x,t,s)=\beta_{n,k}(x,s)$ is defined by \eqref{beta sol}
with $p_{n}=p_n(x;s)$.
\end{lemma}
\begin{proof}
Let $\beta_{n,k}(x,s)$ be defined by \eqref{beta def} with initial
values $\beta_{n,0}(x,s)=p_n(x;s)$. Combining Lemma \ref{L1} with
Proposition \ref{PP3} we see that the initial condition \eqref{H
ini} holds. Therefore, to conclude the proof we only need to verify
the following.
\begin{claim} The matrix $\{\beta_{n,k}(x,s): 0\leq k\leq n\}$ as defined by \eqref{beta def} with $p_{n}=p_n(x;s)$
satisfies equation \eqref{H}.
\end{claim}
A straightforward   computational proof goes as follows. Equation
\eqref{H} is
\begin{multline}
  \label{H2}
  a_n \beta_{n+1,k}(x,s)+b_n(t)\beta_{n,k}(x,s)+c_n(t)\beta_{n-1,k}(x,s)
  =A_{k-1}(x,s)\beta_{n,k-1}(x,s)\\+B_k(x,t,s)\beta_{n,k}(x,s)+C_{k+1}(t,s)\beta_{n,k+1}(x,s).
\end{multline}
In view of \eqref{beta def}, and using the explicit form \eqref{bb},
\eqref{cc} we verify that the coefficients at variable $t$ of this
equation cancel out. Therefore, in \eqref{H2} without loss of
generality we may take $t=0$. We now write $\beta_{n,k}(x,s)$  as
$\sum_{j=0}^k\gamma_{n,k,j}p_{n-j}(x;s)$ where according to
\eqref{beta sol}, we have
\begin{equation}\label{gamma}
\gamma_{n,k,j}=q^{(k-j)(k-j-1)/2}\frac{[n]_q!}{[n-k]_q![j]_q![k-j]_q!}.
\end{equation}
(We will also use the conventions that  $\gamma_{n,k,j}=0$ unless
$0\leq k\leq n$ and $0\leq j \leq k$.) Then \eqref{H2} is equivalent
to a number of identities that arise from comparing the coefficients
at $p_{n-j}(x;s)$. Here we use \eqref{p-rec} to rewrite the terms
$A_k(x,s)p_{n-j}(x;s)$ and $B_k(x,0,s)p_{n-j}(x;s)$ as the linear
combinations of the polynomials $\{p_r(x;s)\}$. We get
\begin{multline}
  \label{H3}
a_n\gamma_{n+1,k,j+1}+b_n(0)\gamma_{n,k,j}+c_n(0)\gamma_{n-1,k,j-1}=
a_{k-1}\gamma_{n,k-1,j}\\-s\eta
q^{k-1}[k-1]_q\gamma_{n,k-1,j}+q^{k-1}a_{n-j-1}\gamma_{n,k-1,j+1}
+q^{k-1}b_{n-j}(s)\gamma_{n,k-1,j}\\
+q^{k-1}c_{n+1-j}(s)\gamma_{n,k-1,j-1}+b_k(0)\gamma_{n,k,j}-(1+q)q^{k-1}s\eta[k]_q\gamma_{n,k,j}\\
+q^k a_{n-j-1}(s)\gamma_{n,k,j+1}+q^kb_{n-j}(s)\gamma_{n,k,j}\\
+q^k c_{n+1-j}(s)\gamma_{n,k,j-1}+c_{k+1}(0)\gamma_{n,k+1,j}-s\eta
q^k [k+1]_q\gamma_{n,k+1,j}.
\end{multline}
Using the identities
\begin{eqnarray*}
  \frac{\gamma_{n,k-1,j}}{\gamma_{n,k,j}}=q^{j-k+1}\frac{[k-j]_q}{[n-k+1]_q},&&\;
\frac{\gamma_{n,k-1,j-1}}{\gamma_{n,k,j}}=\frac{[j]_q}{[n-k+1]_q},\\
\frac{\gamma_{n,k+1,j}}{\gamma_{n,k,j}}=q^{k-j}\frac{[n-k]_q}{[k+1-j]_q},&&\;
\frac{\gamma_{n+1,k,j+1}}{\gamma_{n,k,j}}=q^{j-k+1}\frac{[n+1]_q[k-j]_q}{[n-k+1]_q[j+1]_q},\\
\frac{\gamma_{n-1,k,j-1}}{\gamma_{n,k,j}}=q^{k-j}\frac{[n-k]_q[j]_q}{[n]_q[k+1-j]_q},&&\;
\frac{\gamma_{n,k,j+1}}{\gamma_{n,k,j}}=q^{j-k+1}\frac{[k-j]_q}{[j+1]_q},\\
\frac{\gamma_{n,k,j-1}}{\gamma_{n,k,j}}=q^{k-j}\frac{[j]_q}{[k+1-j]_q},&&\;
\frac{\gamma_{n,k-1,j+1}}{\gamma_{n,k,j}}=q^{2j-2k+3}\frac{[k-j]_q[k-j-1]_q}{[n-k+1]_q[j+1]_q},
\end{eqnarray*} equation \eqref{H3} reduces to the following two identities
between $q$-numbers. The first identity comes from comparing the
coefficients at $s$,
\begin{multline}\label{HHs}
0=-q^{j}\frac{[k-1]_q[k-j]_q}{[n-k+1]_q}+q^{j}\frac{[n-j]_q[k-j]_q}{[n-k+1]_q}+q^{k-1}\frac{[n+1-j]_q[j]_q}{[n-k+1]_q}\\
-(1+q)q^{k-1}[k]_q+q^k[n-j]_q+q^{2k-j}\frac{[n+1-j]_q[j]_q}{[k+1-j]_q}-q^{2k-j}\frac{[k+1]_q[n-k]_q}{[k+1-j]_q}.
\end{multline}
The second identity arises from comparing the coefficients free of
$s$,
\begin{multline}
  \label{HH0}
  a_nq^{j-k+1}\frac{[n+1]_q[k-j]_q}{[n+1-k]_q[j+1]_q}+b_n(0)-b_k(0)+c_n(0)q^{k-j}\frac{[n-k]_q[j]_q}{[n]_q[k+1-j]_q}\\
  =a_{k-1}q^{j-k+1}\frac{[k-j]_q}{[n+1-k]_q}+a_{n-j-1}q^{2j-k+2}\frac{[k-j]_q[k-j-1]_q}{[n+1-k]_q[j+1]_q}\\
  +b_{n-j}(0)q^j\frac{[k-j]_q}{[n+1-k]_q}+c_{n+1-j}(0)q^{k-1}\frac{[j]_q}{[n+1-k]_q}+a_{n-j-1}q^{j+1}\frac{[k-j]_q}{[j+1]_q}\\
  +b_{n-j}(0)q^k+c_{n+1-j}(0)q^{2k-j}\frac{[j]_q}{[k+1-j]_q}+c_{k+1}(0)q^{k-j}\frac{[n-k]_q}{[k+1-j]_q}.
\end{multline}
Identities \eqref{HHs} and \eqref{HH0} are in the form suitable for
computer-assisted verification. We used Mathematica to confirm their
validity.
\end{proof}

\subsection{Monic polynomials}
Let $\overline{Q}_n(y;x,t,s)$ denote the monic version of polynomials
$Q_n$; these polynomials satisfy the recurrence
\begin{multline}
\label{Q rec monic}
y\overline{Q}_n(y;x,t,s)=\overline{Q}_{n+1}(y;x,t,s)+B_n(x,t,s)\overline{Q}_n(y;x,t,s)\\+A_{n-1}(x,s)C_n(t,s)\overline{Q}_{n-1}(y;x,t,s),\;
n\geq 0,
\end{multline}
with the usual initial conditions $\overline{Q}_{-1}=0$, $\overline{Q}_{0}=1$.
Here $A_n,B_n,C_n$ are defined by \eqref{AA}, \eqref{BB}, and
\eqref{CC}, respectively. Let $\overline{p}_n(y;t)=\overline{Q}_n(y;0,t,0)$ be
the monic version of polynomials $p_n$. The monic polynomials
$\{\overline{Q}_n\}$ are well defined for all $x,s$, leading to the
following  version of Lemma \ref{L5}.
\begin{corollary}\label{C bar p2bar Q} 
Suppose $0\leq q\leq 1$, $\eta>0$, $\theta\geq 0$, $\tau\geq 0$.
 For all $s,t>0$, $x,y\in\sR$ we have
\begin{equation}
  \label{p-bar 2 Q-bar}
  \overline{p}_n(y;t)=\sum_{k=0}^n\overline{\beta}_{n,k}(x,s)\overline{Q}_k(y;x,t,s),
  \; n\geq 0.
\end{equation}
\end{corollary}
\begin{proof} Suppose $x\not\in \{x_0(s),x_1(s),\dots\}$. Then
polynomials $Q_n(y;x,t,s)$ are well defined and from Lemma \ref{L5}
we know that \eqref{p2Q} holds with
$\beta_{n,k}(x,t,s)=\beta_{n,k}(x,0,s)$ which do not depend on $t$.

It is well known that the monic polynomials $\overline{Q}_n$ can be
written as
\begin{equation}
  \label{Q-bar2Q}
  \overline{Q}_n(y;x,t,s)=\frac{Q_n(y;x,t,s)}{\prod_{j=0}^{n-1}A_j(x,s)}.
\end{equation}
Therefore for such $x\not\in \{x_0(s),x_1(s),\dots\}$, we get \eqref{p-bar 2 Q-bar} with
\begin{equation}
  \label{bar beta}
  \overline{\beta}_{n,k}(x,s)=\frac{\beta_{n,k}(x,0,s)\prod_{j=0}^{k-1}A_j(x,s)}{\prod_{j=0}^{n-1}a_j}.
\end{equation}

We now extend this relation to all $x$.  
From \eqref{beta def} and \eqref{bar beta} we see that  \eqref{p-bar
2 Q-bar} is a relation between the polynomials in variable $x$ and
holds on an infinite set of $x$. Therefore, it extends to all
$x\in\sR$.
\end{proof}
\subsection{Uniqueness}
It turns out that polynomials $\{Q_n\}$ can be interpreted as
 polynomials $\{p_n\}$ with modified parameters.
\begin{lemma}\label{L4} Suppose $x\not\in\{x_0(s),x_1(s),\dots\}$.
With
\begin{equation}\label{prim}
 \theta'=\frac{\theta-s\eta-(1-q)x}{1+\eta x},\;
\tau'=\frac{\tau+(1-q)s}{1+\eta x},
\end{equation}
we have
\begin{equation}\label{p-Q}
  Q_n^{(\theta,\tau)}(y;x,t,s)=p_n^{(\theta',\tau')}\left(\frac{y-x}{1+\eta x};\frac{t-s}{1+\eta
x}\right).
\end{equation}
\end{lemma}
\begin{proof}
Write \eqref{Q rec} as
\begin{multline*}
  \frac{y-x}{1+\eta x}Q_n(y)=\frac{A_n(x,s)}{1+\eta x}
  Q_{n+1}(y)+\frac{B_n(x,t,s)-x}{1+\eta x}Q_n(y)\\+\frac{C_n(t,s)}{1+\eta
  x}Q_{n-1}(y).
\end{multline*}
Consider polynomials $r_n(y')$ such that $r_{-1}=0$, $r_0=1$, and
\begin{multline*}
  y'r_n(y')=\frac{A_n(x,s)}{1+\eta x}
  r_{n+1}(y')+\frac{B_n(x,t,s)-x}{1+\eta x}r_n(y')+\frac{C_n(t,s)}{1+\eta
  x}r_{n-1}(y').
\end{multline*}
Since $r_{-1}=Q_{-1}$ and $r_0=Q_0$, setting $y'=\frac{y-x}{1+\eta
x}$ we have
$$ r_n(y')=Q_n(y),\; n\geq 1.
$$
Using \eqref{prim} we get
\begin{multline*}
  y'r_n(y')=\left(\eta^{-1}+\theta'[n]_q+\tau'\eta[n]_q^2\right)
  r_{n+1}(y')\\+\left(\eta t'+\theta'+([n]_q+[n-1]_q)\tau'\eta\right)[n]_qr_n(y')
  \\+\eta\left(t'+\tau'[n-1]_q\right)[n]_qr_{n-1}(y').
\end{multline*}
This means that polynomials $r_n(y')$ satisfy the same recurrence as
polynomials $p_n(y';t')$ with parameters $\theta'$ and $\tau'$. Thus
\eqref{p-Q} follows.
\end{proof}

Polynomials $\{Q_n\}$ are just a reparametrized version of
polynomials $p_n$, see Proposition \ref{PP3}, so their orthogonality
measure $\nu_{x,t,s}$ is also determined by moments. Since the
orthogonality measure of polynomials $\overline{Q}_n$ may differ only for
$x\in\{x_0(s),x_1(s),\dots\}$ in which case it has finite support,
we get the following.
\begin{corollary}\label{C1}
For all $x$ such that \eqref{Prod A} holds, the orthogonality
measure $\nu_{x,t,s}$ of polynomials $\overline{Q}(\cdot;x,t,s)$ is
unique.
\end{corollary}

\subsection{Proof of Theorem \protect{\ref{T1}}}\label{Proof of T1}
\begin{proof}[Proof of Theorem \ref{T1}]
Replacing $x$ by $-x$ in \eqref{p-monic}, changes $\eta,\theta$ to
$-\eta,-\theta$. So without loss of generality we may assume
$\eta\geq 0$. Furthermore, the case $\eta=0$ is known from
\cite{Bryc-Wesolowski-03}, so we only consider $\eta>0$.

 Let
$\nu_{x,t,s}(dy)$ be the orthogonality measure of polynomials
$\overline{Q}_n(y;x,t,s)$, see \eqref{Q rec monic}. By Proposition
\ref{P+},  measure $\nu_{x,t,s}(dy)$ is well defined for all
$0<s<t$, $x\in \supp(\mu_s)$.

 Corollary \ref{C bar
p2bar Q} implies that
$$\overline{\beta}_{n,k}(x,s)=\frac{\int
\overline{p}_n(y;t)\overline{Q}_k(y;x,t,s)\nu_{x,t,s}(dy)}{\|\overline{Q}_k\|_2^2}.$$
Since $\overline{p}_n(x;s)=\overline{\beta}_{n,0}(x,s)$, using the above with
$k=0$ we see that projection formula \eqref{proj main} holds for all
$x\in \supp(\mu_s)$.

Projection formula \eqref{proj main} determines the moments of
$\nu_x$. By Corollary \ref{C1}, this determines $\nu_x$
uniquely.
\end{proof}
\subsection*{Acknowledgement} The first-named author (WB) thanks
Mourad Ismail for several helpful conversations. The authors are grateful to the anonymous referees
for {\em Constructive Approximation} who helped us to avoid an embarrassing blunder.


\end{document}